\documentclass[12pt,draft]{amsart}

\makeatletter

\theoremstyle{plain}
\newtheorem{thm}{Theorem}[section]
\numberwithin{equation}{section} 
\numberwithin{figure}{section} 
\theoremstyle{plain}
\newtheorem*{thm*}{Theorem}
\theoremstyle{plain}
\newtheorem{cor}[thm]{Corollary} 
\theoremstyle{plain}
\newtheorem{lem}[thm]{Lemma} 
\theoremstyle{plain}
\newtheorem{prop}[thm]{Proposition} 
\theoremstyle{definition}
\newtheorem{defn}[thm]{Definition}
\theoremstyle{remark}

\theoremstyle{remark}

\theoremstyle{remark}

\theoremstyle{remark}

\theoremstyle{definition}

\theoremstyle{remark}
  \newtheorem*{acknowledgement*}{Acknowledgement}

\theoremstyle{plain}

\theoremstyle{plain}
\theoremstyle{plain}
\theoremstyle{plain}
\theoremstyle{definition}

\theoremstyle{remark}

\theoremstyle{remark}

\theoremstyle{remark}

\theoremstyle{plain}

\newcommand{\mtp}{\mathop{\otimes}_{\max}}

\newcommand{\C}{{\mathbb C}}

\newcommand{\M}{{\mathbb M}}

\newcommand{\B}{{\mathbb B}}

\newcommand{\p}{\varphi}
\newcommand{\hh}{{\mathcal H}}
\newcommand{\hk}{{\mathcal K}}

\newcommand{\tr}{\mathop{\mathrm{tr}}\nolimits}
\newcommand{\Tr}{\mathop{\mathrm{Tr}}\nolimits}

\usepackage{amssymb} \usepackage{amscd}

\makeatother

\begin{document}

\title[Property T]{Kazhdan's Property T and C$^*$-algebras}

\author{Nathanial P. Brown}

\address{Department of Mathematics, Penn State University, State
College, PA 16802}

\email{nbrown@math.psu.edu}

\thanks{Partially supported by DMS-0244807.}

\begin{abstract}
Kazhdan's property T has recently been imported to the C$^*$-world by
Bekka.  Our objective is to extend a well known fact to this realm;
we show that a nuclear C$^*$-algebra with property T is finite
dimensional (for all intents and purposes). Though the result is not
surprising, the proof is a bit more complicated than the group case.
\end{abstract}

\maketitle

\section{Introduction}

Kazhdan's revolutionary concept of property T has recently been
translated into C$^*$-language in \cite{bekka}.  One of the questions
raised by Bekka's paper is whether or not one can generalize to the
C$^*$-context the classical fact that a discrete group which is both
amenable and has property T must be finite (cf.\ \cite[Proposition
11]{bekka}). Unfortunately, the C$^*$-situation is not quite as
simple, but a satisfactory result can be obtained.

\begin{thm*} Let $A$ be a unital C$^*$-algebra which is both nuclear
and has property T.  Then $A = B \oplus C$ where $B$ is finite
dimensional and $C$ has no tracial states.
\end{thm*}

The irritating $C$-summand can't be avoided; if $B$ is {\em any}
C$^*$-algebra with property T and $C$ is {\em any} algebra without
tracial states then $B \oplus C$ also has property T.  Hence any
finite dimensional C$^*$-algebra plus a Cuntz algebra (for example)
will have property T and be nuclear.  On the other hand, the theorem
above does imply that if $A$ is nuclear, has property T and has a
{\em faithful} trace then it must be finite dimensional -- this is an
honest generalization of the discrete group case since reduced group
C$^*$-algebras always have a faithful trace.  More generally, every
stably finite, nuclear algebra with property T is finite dimensional
since Haagerup has shown that every unital stably finite exact
C$^*$-algebra must have a tracial state (cf.\
\cite{haagerup-thorbjornsen}). Since these are the main cases of
interest, it seems fair to say {\em ``property T plus amenability
implies finite dimensional} (more or less)."

Perhaps the more interesting thing, however, is the proof. In the
case of a discrete group it is trivial: amenability implies the left
regular representation has almost invariant vectors; rigidity then
provides a fixed vector; but, only finite groups have fixed vectors
in the left regular representation.

Unfortunately, we have been unable to find a simple argument for the
general case, hence the circuitous route taken here.  Our approach
requires generalizing Kazhdan projections, the theory of amenable
traces and even the deep fact that nuclearity passes to quotients.

Finally, we express our gratitude to the reviewer for pointing out an 
error in our original manuscript; their careful reading greatly 
improved the truth of this work! 

\section{Definitions and Notation}

We make the blanket assumption that {\em all C$^*$-algebras are
unital and separable} unless otherwise noted or obviously false
(e.g.\ $\B(\hh)$, the bounded operators on a (separable) Hilbert
space $\hh$, won't be norm separable).

Inspired by the von Neumann version (see, for example,
\cite{connes-jones}), Bekka defines property T for C$^*$-algebras in
terms of bimodules. If $A$ is a C$^*$-algebra and $\hh$ is a Hilbert
space equipped with commuting actions of $A$ and its opposite algebra
$A^{op}$ then we say $\hh$ is an $A$-$A$ bimodule.  (Another way of
saying this is that there exists a $\ast$-representation $\pi\colon
A\mtp A^{op} \to \B(\hh)$.) As is standard, we denote the action by
$\xi \mapsto a\xi b$, $\xi \in \hh$, $a \in A$, $b \in A^{op}$. (That
is, $a\xi b = \pi(a\otimes b)\xi$.)

\begin{defn} A C$^*$-algebra $A$ has property T if every bimodule
with almost central vectors has a central vector; i.e.\ if $\hh$
is a bimodule and there exist unit vectors $\xi_n \in \hh$ such
that $\|a\xi_n - \xi_n a\| \to 0$ for all $a \in A$ then there
exists a unit vector $\xi \in \hh$ such that $a\xi = \xi a$ for
all $a \in A$.
\end{defn}

An important example of a bimodule is gotten by starting with an
embedding $A \subset \B(\hk)$ and letting $\hh = \mathrm{HS}(\hk)$
be the Hilbert-Schmidt operators on $\hk$.  The commuting actions
of $A$ and $A^{op}$ are given by multiplication on the left and
right: $T \mapsto aTb$ for all $T \in \mathrm{HS}(\hk)$, $a \in A$
and $b \in A^{op}$ (canonically identified, as normed involutive
linear spaces, with $A$).

Throughout this note we will use $\Tr$ to denote the canonical
(unbounded) trace on $\B(\hh)$ and, if $\hh$ happens to be finite
dimensional, $\tr$ will be the unique tracial state on $\B(\hh)$.

\section{Kazhdan Projections}

It is known that if $\Gamma$ is a discrete group with property T
then all the Kazhdan projections -- the central covers in the
double dual $C^*(\Gamma)^{**}$ coming from finite dimensional
irreducible representations -- actually live in $C^*(\Gamma)$. We
extend this fact to the general C$^*$-context.

Recall that an {\em intertwiner} of two $*$-representations
$\pi\colon A\to \B(\hh)$ and $\sigma\colon A\to \B(\hk)$ is a
bounded linear operator $T\colon \hh\to \hk$ such that $T\pi(a) =
\sigma(a)T$ for all $a \in A$. We will need {\em Schur's Lemma}.

\begin{lem}
\label{lem:Schur} If two representations $\pi$ and $\sigma$ have a
nonzero intertwiner and $\pi$ is irreducible then $\pi$ is
unitarily equivalent to a subrepresentation of $\sigma$.
\end{lem}

The proof is simple, well-known and will be omitted -- the main
point is that irreducibility of $\pi$ forces an intertwiner to be
a scalar multiple of an isometry.

Property T groups are often defined as follows: if a unitary
representation weakly contains the trivial representation then it
must honestly contain it.  Here is the generalization to our
context.

\begin{prop}
\label{prop:fdweakcontain} Assume $A$ has property T, $\pi\colon A
\to \M_n({\mathbb C})$ is an irreducible representation and
$\sigma\colon A \to \B(\hk)$ is any representation which weakly
contains $\pi$.\footnote{This means there exist isometries $V_k\colon
{\mathbb C}^n \to \hk$ such that $\|\sigma(a)V_k - V_k \pi(a)\| \to
0$ for all $a \in A$. This is equivalent to saying
$\|V_k^*\sigma(a)V_k - \pi(a)\| \to 0$ for all $a \in A$, which
explains the terminology `weak containment'.}  Then $\pi$ is
unitarily equivalent to a subrepresentation of $\sigma$.
\end{prop}

\begin{proof}
Let $\mathrm{HS}({\mathbb C}^n,\hk)$ denote the Hilbert-Schmidt
operators from ${\mathbb C}^n$ to $\hk$.  We make this space into an
$A$-$A$ bimodule by multiplication on the left and right -- i.e.\
$aTb = \sigma(a)T\pi(b)$ for all $T \in \mathrm{HS}({\mathbb
C}^n,\hk), a \in A, b \in A^{op}$.  Since a nonzero central vector
would evidently be an intertwiner of $\pi$ and $\sigma$, it suffices
(by property T and Schur's lemma) to show the existence of an
asymptotically central sequence in $\mathrm{HS}({\mathbb C}^n,\hk)$.

If $V_k\colon {\mathbb C}^n \to \hk$ are isometries such that
$\|\sigma(a)V_k - V_k \pi(a)\| \to 0$ then a routine calculation
shows that the unit vectors $\frac{1}{\sqrt{n}}V_k \in
\mathrm{HS}({\mathbb C}^n,\hk)$ have the property that $\| a
\frac{1}{\sqrt{n}}V_k - \frac{1}{\sqrt{n}}V_k a \|_{\mathrm{HS}}
\to 0$ for all $a \in A$.
\end{proof}

Recall that if $\pi\colon A \to \B(\hh)$ is a representation then
there is a central projection $c(\pi) \in \mathcal{Z}(A^{**})$ in
the double dual of $A$ with the property that $c(\pi)A^{**} \cong
\pi(A)^{\prime\prime}$ (among other things -- see \cite{pedersen}
for more).  For property T groups and finite dimensional
representations these projections are often called something else.

\begin{defn} If $A$ has property T and $\pi\colon A \to \M_n(\C)$
is irreducible then the central cover $c(\pi)$ is also known as
the {\em Kazhdan projection} associated to $\pi$.
\end{defn}

One of the remarkable consequences of property T is that all
Kazhdan projections actually live in (the center of) $A$ (not just
$A^{**}$). (Compare with the fact that $C^*(\mathbb{F}_2)$ has
tons of finite dimensional representations, yet no nontrivial
projections!)

\begin{thm} Assume $A$ has Kazhdan's property T.  Then
for each finite dimensional irreducible representation $\pi\colon
A \to \M_n({\mathbb C})$, the Kazhdan projection $c(\pi)$ actually
lives in $A$.
\end{thm}

\begin{proof} Let $\sigma\colon A \to \B(\hh)$ be a representation
with the following three properties: (1) $\sigma(A)$ contains no
nonzero compact operators, (2) $\pi \oplus \sigma\colon A \to
\B({\mathbb C}^n \oplus \hh)$ is faithful and (3) $\sigma$
contains no subrepresentation which is unitarily equivalent to
$\pi$. For example, one can start with a faithful representation
of the algebra $(1 - c(\pi))A \subset A^{**}$ and inflate, if
necessary, to arrange (1). (Standard theory of central covers
shows that such a $\sigma$ has no subrepresentation unitarily
equivalent to $\pi$.)

Notice that such a representation $\sigma$ can't possibly be
faithful -- if it were then Voiculescu's Theorem (cf.\
\cite{davidson}) would imply that $\sigma$ is approximately
unitarily equivalent to $\sigma \oplus \pi$. In other words,
$\sigma$ weakly contains $\pi$ and thus, by Proposition
\ref{prop:fdweakcontain}, actually contains $\pi$. This
contradicts our assumption (3) and so $\sigma$ can't be faithful.

Thus $J = ker(\sigma)$ is a nontrivial ideal in $A$.  But
assumption (2) implies that $\pi|_J$ must be faithful; hence, $J$
is finite dimensional and has a unit $p$ which is necessarily a
central projection in $A$.  We will show $p = c(\pi)$.

Since $\pi$ is irreducible, $\pi(p) = 1$ and $J \cong \M_n(\C)$.
Hence we may identify the representations $A \to pA$ and $\pi$.
Thus they have the same central covers -- i.e.\ $p = c(\pi)$ as
desired.
\end{proof}

Here are a couple of consequences.

\begin{cor} If $A$ has property T then it has at most
countably many non-equivalent finite dimensional representations.
\end{cor}

\begin{proof} A separable C$^*$-algebra has at most countably many
orthogonal projections (since it can be represented on a separable
Hilbert space).
\end{proof}

\begin{cor}
\label{cor:nofdquotient} Assume $A$ has property T and let $J
\triangleleft A$ be the ideal generated by all of the Kazhdan
projections.  Then $A/J$ has no finite dimensional
representations.
\end{cor}

\begin{proof}  Any nonzero, finite dimensional representation of
$A/J$ would produce a Kazhdan projection in $A$ which wasn't in the
kernel of the quotient map $A \to A/J$. \end{proof}

\section{Amenable Traces}

The notion of amenability for traces has a reasonably long history,
with important contributions from several authors (see \cite{brown}
for history and references). In this section we adapt one of
Kirchberg's contributions (cf.\ \cite{kirchberg}).

\begin{defn}
\label{defn:amenabletrace} Let $A \subset \B(\hh)$ be a concretely
represented unital C$^*$-algebra.  A tracial state $\tau$ on $A$ is
called {\em amenable} if there exists a state $\p$ on $\B(\hh)$ such
that (1) $\p|_A = \tau$ and (2) $\p(uTu^*) = \p(T)$ for every unitary
$u \in A$ and $T \in \B(\hh)$.
\end{defn}

It is a remarkable fact (due to Connes and Kirchberg) that this
notion can be recast in terms of approximation by finite dimensional
completely positive maps. See \cite{brown} or \cite{kirchberg} for a
proof of the theorem and \cite{paulsen:book} for more on completely
positive maps.

\begin{thm}
\label{thm:amenabletraces} Let $\tau$ be a tracial state on $A$. Then
$\tau$ is amenable if and only if there exist unital completely
positive maps $\p_n \colon A \to \M_{k(n)} (\C)$ such that $\| \p_n
(ab) - \p_n (a) \p_n (b) \|_{2,\mathrm{tr}} \to 0$, where
$\|x\|_{2,\mathrm{tr}} = \sqrt{\mathrm{tr}(x^*x)}$, and $\tau (a) =
\lim_{n \to \infty} \mathrm{tr} \circ \p_n (a)$, for all $a,b \in A$.
\end{thm}

With this approximation property in hand, the following proposition
is straightforward.

\begin{prop}
\label{prop:fdquotient} Let $A$ be a C$^*$-algebra with property T.
Then $A$ has an amenable trace if and only if $A$ has a nonzero
finite dimensional quotient.
\end{prop}

\begin{proof} Evidently a finite dimensional quotient yields an
amenable trace (since every trace on a finite dimensional
C$^*$-algebra is amenable).  Hence we assume $A$ has property T and
an amenable trace $\tau$.

Let $\p_n \colon A \to \M_{k(n)} (\C)$ be unital completely positive
maps as in Theorem \ref{thm:amenabletraces}. Invoking Stinespring's
Theorem, we can find representations $\rho_n\colon A \to \B(\hh_n)$
and finite rank projections $P_n \in \B(\hh_n)$ such that $\p_n$ can
be identified with $x \mapsto P_n\rho_n(x)P_n$. Let $$\rho =
\bigoplus_{n \in \mathbb{N}} \rho_n \colon A \to \B(\bigoplus_{n \in
\mathbb{N}} \hh_n)$$ and regard $\mathrm{HS}(\bigoplus_{n \in
\mathbb{N}} \hh_n)$ as an $A$-$A$ bimodule via $T \mapsto
\rho(a)T\rho(b)$.

The identity $P_n\rho_n(a) - \rho_n(a)P_n = P_n\rho_n(a)(1-P_n) -
(1-P_n)\rho_n(a)P_n$ together with an unenlightening calculation
shows that
$$\frac{\|P_n\rho_n(a) - \rho_n(a)P_n\|_2}{\|P_n\|_2} =
\bigg(\mathrm{tr}(\p_n(aa^*) - \p_n(a)\p_n(a^*)) +
\mathrm{tr}(\p_n(a^*a) - \p_n(a^*)\p(a)) \bigg)^{\frac{1}{2}},$$
where $\|\cdot \|_2$ denotes the Hilbert-Schmidt norm induced by
$\langle S,T\rangle = \Tr(T^*S)$. By the Cauchy-Schwarz inequality we
have $$\mathrm{tr}(\p_n(aa^*) - \p_n(a)\p_n(a^*)) \leq \| \p_n (aa^*)
- \p_n (a) \p_n (a^*) \|_{2,\mathrm{tr}}$$ and hence
$$\frac{\|P_n\rho_n(a) - \rho_n(a)P_n\|_2}{\|P_n\|_2} \to 0,$$ for
every $a \in A$.  That is, $\mathrm{HS}(\bigoplus_{n \in \mathbb{N}}
\hh_n)$ has a sequence of asymptotically central unit vectors
(namely, $\frac{1}{\|P_n\|_2}P_n$). Hence property T gives us a
nonzero central vector $T \in \mathrm{HS}(\bigoplus_{n \in
\mathbb{N}} \hh_n)$ -- i.e.\ a nonzero compact operator in the
commutant of $\rho(A)$. Since $T$'s spectral projections must also
live in $\rho(A)^{\prime}$ we get a finite rank projection in the
commutant. Thus $\rho(A)$, and hence $A$, has a finite dimensional
quotient.
\end{proof}

Since property T evidently passes to quotients, the following
corollary is a consequence of the previous result and Corollary
\ref{cor:nofdquotient}.

\begin{cor}
\label{cor:noamenabletrace} Assume $A$ has property T and let $J
\triangleleft A$ be the ideal generated by all the Kazhdan
projections. Then $A/J$ has no amenable traces.
\end{cor}

\section{Nuclearity and Property T}

We now have the ingredients necessary for the mundane fact that
nuclear C$^*$-algebras with property T are finite dimensional plus
something traceless.

\begin{thm} Assume $A$ is nuclear and has property T.  Then $A = B
\oplus C$ where $B$ is finite dimensional and $C$ admits no
tracial states.
\end{thm}

\begin{proof} Let $J$ be the ideal generated by the Kazhdan
projections in $A$. Note that $A/J$ has property T and is nuclear.
Since every trace on a nuclear C$^*$-algebra is amenable (cf.\
\cite{brown}), Corollary \ref{cor:noamenabletrace} implies that $A/J$
is traceless.  Thus it suffices to show $B = J$ is finite dimensional
(since it would then have to be a direct summand and $C= A/J$ is
traceless).

So, how to see that $J$ is finite dimensional?  Well, if it
weren't then we could find integers $k(n)$ such that $$J =
\bigoplus_{n=1}^{\infty} \M_{k(n)}(\C),$$ the $c_0$-direct sum
(sequences tending to zero in norm) and thus the multiplier
algebra of $J$ is equal to the algebra of bounded sequences. Hence
there is a unital $*$-homomorphism $$A/J
\to\frac{\prod_{n=1}^{\infty}
\M_{k(n)}(\C)}{\bigoplus_{n=1}^{\infty}\M_{k(n)}(\C)}.$$  However,
the Corona algebra on the right is easily seen to have lots of
tracial states and so we deduce that $A/J$ has a tracial state.
But this is silly, as observed in the preceding paragraph.
\end{proof}

\bibliographystyle{amsplain}

\begin{thebibliography}{10}

\bibitem{bekka} B. Bekka, \emph{A property (T) for C$^*$-algebras},
preprint 2005.

\bibitem{brown} N.P. Brown, \emph{Invariant means and finite representation
theory for C$^*$-algebras}, Mem.\ Amer.\ Math.\ Soc.\ (to appear).

\bibitem{connes-jones} A.~Connes and V.~Jones, \emph{Property T for
von Neumann algebras}, Bull.\ London Math.\ Soc.\ \textbf{17}
(1985), 57--62.

\bibitem{davidson} K.R. Davidson, \emph{$C\sp *$-algebras by example},
Fields Institute Monographs, 6. American Mathematical Society, Providence, RI, 1996.

\bibitem{haagerup-thorbjornsen} U. Haagerup and S. Thorbj{\o}rnsen,
\emph{Random matrices and K-theory for exact C$^*$-algebras}, Doc.\ Math.\
\textbf{4} (1999), 341--450.

\bibitem{kirchberg} E. Kirchberg, \emph{Discrete groups with Kazhdan's property
${\rm T}$ and factorization property are residually finite},  Math.\ Ann.\  \textbf{299}  (1994),
551--563.

\bibitem{paulsen:book} V.~Paulsen, \emph{Completely bounded maps and
operator algebras}, Cambridge Studies in Advanced Mathematics 78,
Cambridge University Press, Cambridge, 2002.

\bibitem{pedersen} G.K. Pedersen, \emph{$C\sp{*} $-algebras and their automorphism groups},
London Mathematical Society Monographs, 14. Academic Press, Inc., London-New York, 1979.

\end{thebibliography}
\providecommand{\bysame}{\leavevmode\hbox to3em{\hrulefill}\thinspace}

\end{document}